\documentclass[number,citesort,MSNbibl,dvips]{arxbj}

\aid{0}
\volume{18}
\issue{2}
\pubyear{2012}
\firstpage{635}
\lastpage{643}
\doi{10.3150/10-BEJ336}

\makeatletter

\newtheorem{theorem}{Theorem}
\newremark{remark}{Remark}
\newtheorem{proposition}{Proposition}
\newtheorem{corollary}{Corollary}

\makeatother

\begin{document}
\begin{frontmatter}

\title{A unified minimax result for restricted parameter spaces}
\runtitle{A unified minimax result}

\begin{aug}
\author[a]{\fnms{\'{E}ric} \snm{Marchand}\corref{}\thanksref{a}\ead[label=e1]{eric.marchand@usherbrooke.ca}} \and
\author[b]{\fnms{William E.} \snm{Strawderman}\thanksref{b}\ead[label=e2]{straw@stat.rutgers.edu}}
\runauthor{\'{E}. Marchand and W.E. Strawderman}
\address[a]{D\'{e}partement de math\'{e}matiques, Universit\'{e} de
Sherbrooke, Qu\'ebec, Canada J1K 2R1.\\ \printead{e1}}
\address[b]{Department of Statistics and Biostatistics, Rutgers
University, 561 Hill
Center, Busch Campus, Piscataway, NJ 08854-8019, USA. \printead{e2}}
\end{aug}

\received{\smonth{4} \syear{2010}}
\revised{\smonth{10} \syear{2010}}

%
\begin{abstract}
We provide a development that unifies, simplifies and
extends considerably a number of minimax results in the restricted
parameter space literature. Various applications follow, such as
that of estimating location or scale parameters under a lower (or
upper) bound restriction, location parameter vectors restricted to
a polyhedral cone, scale parameters subject to restricted ratios
or products, linear combinations of restricted location
parameters, location parameters bounded to an interval with
unknown scale, quantiles for location-scale families with
parametric restrictions and restricted covariance matrices.
\end{abstract}

%
\begin{keyword}
\kwd{covariance matrices}
\kwd{linear combinations}
\kwd{location parameters}
\kwd{minimax}
\kwd{polyhedral cones}
\kwd{quantiles}
\kwd{restricted parameters}
\kwd{scale parameters}
\end{keyword}

\end{frontmatter}

\section{Introduction}\label{sec1}
We provide a development that unifies, simplifies and extends
considerably a number of minimax results in the restricted
parameter space estimation literature. As illustrated with a
series of examples, the unified minimax result has wide
applicability with respect to the nature of the constraint, the
underlying probability model and the loss function utilized.

To further put into context the findings of this paper, consider a
basic situation where $X \sim N(\theta,1)$, with $\theta\geq a$
($a>-\infty$ known), and where $\theta$ is estimated under squared
error loss $(d-\theta)^2$. Katz \cite{Kat61} established that the Bayes
estimator $\delta_U$ with respect to the flat prior on
$(a,\infty)$ dominates the minimum risk equivariant (MRE) estimator
$\delta_0(X)=X$. However,
$\delta_0$ remains a useful benchmark estimator with its constant risk
matching the minimax risk, and with any
improvement, such as $\delta_U$, being necessarily minimax as
well. In a technical sense and roughly speaking, the form (and
unboundedness) of the restricted parameter space $[a,\infty)$
preserves a common structure with the unrestricted parameter space
$\Re$, and the constructions of the least favourable sequence of
priors for both problems are isomorphic, leading to the same
minimax values. In contrast, the restriction to a compact interval
$\theta\in[a,b]$ is quite different and lowers the minimax risk
(see example (C) in Section \ref{sec3}).

Now, the above phenomenon is not only more general (e.g., general
location families with absolutely continuous Lebesgue densities
and strictly convex loss \cite{Far64}), but similar results have
been established in various other situations, beginning with
Blumenthal and Cohen \cite{BluCoh68} in the context of ordered location
parameters. Many other such contributions will be referred to
below, but at this point we refer to the monograph of van Eeden~\cite{van06},
as well as the review paper by Marchand and Strawderman
\cite{MarStr04}, which contain a~substantial amount of material and
references relating to such problems.

In this paper, we provide a unified framework for the
above-mentioned problems, as well as many others either for more
general loss and/or model, or for new situations such as
estimating quantiles or covariance matrices under parametric
restrictions. While results for certain of the problems (e.g.,
(A) and (B)) are not new (although we generalize some to more general
loss functions) and certain others have been studied for squared
error loss, we greatly expand the set of loss
functions for which minimaxity is established (e.g., (C)--(F)), certain of the problems (e.g., (G) and (H) and Remark \ref{rem2}) have
not been extensively studied and thus our problems and results are
mostly new. In Section \ref{sec2}, we formalize the general argument,
relying on the existence of a least favourable sequence
(Proposition \ref{brown}), setting up \textit{conditions} on the
restricted parameter space that facilitate a~correspondence with
the above sequence (Theorem \ref{minimax}) and inferring
(Corollary~\ref{mre}) that a minimax MRE estimator remains minimax with
the introduction of a~restriction on the parameter space under given \textit{conditions}.
Detailed examples follow in Section \ref{sec3}. These include the
estimation of location or scale parameters under a lower (or
upper) bound restriction, location parameter vectors restricted to
a polyhedral cone, scale parameters subject to restricted ratios
or products, linear combinations of restricted location
parameters, location parameters bounded to an interval with
unknown scale, quantiles for location-scale families with
parametric restrictions and covariance matrices with restricted traces
or determinants.

\section{Main result}\label{sec2}

We begin with the following fact concerning minimax problems,
presented as a synthesized version of parts of the Appendix of
\cite{Bro86}, pages 254--268.

\begin{proposition}
\label{brown} Let $R< \infty$ be the minimax value in a problem
with sample space $\textsl{X}$ and parameter space $\Omega$, both
Euclidean. Suppose the probability measures are absolutely
continuous with respect to a $\sigma$-finite measure, that the
loss $L(\theta, \cdot)$ is lower semicontinuous on the
action space and that $L(\theta,a) \to b(\theta)=\sup L(\theta, \cdot)$
as $\|a\| \to\infty$ for all $\theta$. Then there exists a sequence of prior
distributions with finite support and with Bayes risks equal to~$r_n$,
such that $r_n$ approaches $R$ as $n \to\infty$, and there also
exists a minimax
procedure.
\end{proposition}

We make use of a classical framework for invariant statistical
problems. This includes a~group of
transformations $G$ with an invariant family of probability measures $\{
P_{\theta}\dvt \theta\in\Omega\}$,
where $X \sim P_{\theta}$ and $X'=g(X)$ implies $X' \sim P_{\theta'}$
with $\theta'=\bar{g}\theta$ and $\bar{G}=
\{\bar{g}\dvt g \in G \}$ forming a corresponding group of actions on
$\Omega$. As well, for estimating a~parametric function
$\tau(\theta)$ with loss $L$, additional assumptions include the
condition that~$\tau(\bar{g} \theta)$ depends on $\theta$ only through
$\tau(\theta)$,
and that the group action on the decision space~$D$ satisfies the
condition $L(\bar{g} \theta, g^* d) = L(\theta,d)$
for all $\theta,d$ (e.g., \cite{LehCas98}, Section
3.2). With the help of Proposition \ref{brown}, we obtain the following result.
\begin{theorem}
\label{minimax} Let a problem satisfying the conditions of Proposition
\ref{brown} be
invariant under a group $G$ and let $\delta_0(X)$ be minimax for
a full parameter space $\Omega$. Suppose now that the parameter space
is restricted to a subset $\Omega^*$; that there exist
sequences $\bar{g}_{n} \in\bar{G}$ and $B_n \subseteq\Omega^*$,
such that $\bar{g}{_n} B_n \subset\bar{g}_{n+1} B_{n+1}$; and
that $\bigcup_{n} \bar{g}_{n} B_n = \Omega$. Then $\delta_0(X)$
remains minimax in the restricted parameter space problem.
\end{theorem}

\begin{pf}
Let $\pi_n$, $S_n$ and $r_n$ be, respectively,
Proposition \ref{brown}'s sequence of priors, sequence of
corresponding finite supports and Bayes risks, with $r_n \to R$
as $n \to\infty$. Choose $m(n)$ sufficiently large so that $m
\geq m(n)$ implies $\bar{g}{_m} B_m \supset S_n$. As we show
below, the prior distribution with finite support
$S_n^*=\bar{g}_{m(n)}^{-1}(S_n) \subset\Omega^*$, given by
$\pi_n^*(\theta)=\pi_n(\bar{g}_{m(n)}^{-1} \theta),$ has Bayes risk
$r_n^*=r_n$. This implies directly that $\delta_0(X)$ is minimax,
since $r_n^*=r_n \to R$, as $n \to\infty$ by Theorem 5.18 of
\cite{Ber85}. It remains to show that the Bayes risks of $\pi_n$
and $\pi_n^*$ coincide, and a standard argument is as follows. Let
$\delta(X)$ be any estimator. Then, for its risk, we have:
\[
R(\theta, \delta) = E_{\theta} L(\theta, \delta(X)) = E_{\bar{g}\theta}
L(\theta, \delta(g^{-1}X)) = E_{\bar{g}\theta} L(\bar{g} \theta,
g^*\delta(g^{-1}X)) = R(\bar{g} \theta, g^*\delta(g^{-1}X)),
\]
by invariance. It follows that, if we set
$\tilde{\theta}=\bar{g}_{m(n)} \theta$, then
\begin{eqnarray*}
r_n &=&
E^{\theta}[R(\theta, \delta_n(X))] =
E^{\theta}\bigl[R\bigl(\bar{g}_{m(n)}^{-1}\theta, {g_{m(n)}^{*-1}}
\delta_n\bigl(g_{m(n)}X\bigr)\bigr)\bigr] \\
&=& E^{\tilde{\theta}}\bigl[R\bigl(\tilde{\theta},
{g^{*-1}_{m(n)}} \delta_n\bigl(g_{m(n)}(X)\bigr)\bigr)\bigr] = r_n^*,
\end{eqnarray*}
where
$\delta_n(X)$ is the Bayes estimator corresponding to $\pi_n$ and
hence ${g^{*-1}_{m(n)}} \delta_n(g_{m(n)}(X))$ is the Bayes
estimator corresponding to $\pi_n^*$.
\end{pf}

For applications, we will take $B_n$ of Theorem \ref{minimax} to
match $\Omega^*$, but it is potentially more convenient to take
$B_n$ as a sequence of open neighborhoods. Now, since the best
equivariant estimators are often minimax, we deduce the following
widely applicable result.

\begin{corollary}
\label{mre} If an MRE estimator in a given problem satisfying the
conditions of Theorem \ref{minimax} is minimax, then it remains
minimax in the restricted problem provided the restricted
parameter space $\Omega^*$ satisfies the conditions of Theorem
\ref{minimax}.
\end{corollary}

For the sake of clarity, we do not assume that the action space
and the image of the restricted parameter space coincide. Hence,
minimax estimators that can
be derived from Theorem \ref{minimax} or Corollary \ref{mre} are
not forced to take values in $\Omega^*$. The main motivation
resides in the benchmarking (i.e., dominating estimators that take
values in $\Omega^*$ are necessarily minimax) and preservation of
minimaxity (the minimax risks on $\Omega$ and $\Omega^*$ are
equivalent). An important class of further applications of Corollary
\ref{mre} will arise in cases where $\delta_{\mathrm{MRE}}$ is
minimax for the unrestricted problem $\Omega$ and the parameter space
$\Omega^*$ and loss $L(\theta, \cdot)$ are convex, in which case
the projection of $\delta_{\mathrm{MRE}}$ onto $\Omega^*$ will
dominate $\delta_{\mathrm{MRE}}$ and hence be minimax.

\begin{remark}
\label{Hunt-Stein} Notwithstanding the conditions required on the
restricted parameter space~$\Omega^*$, the applicability of
Corollary \ref{mre} hinges on the minimaxity of the best
equivariant estimator, in particular for unrestricted parameter
space versions. As studied and established by several authors, it
turns out that it is frequently the case that a minimax
equivariant rule exists. We refer to \cite{Ber85}, Section 6.7,
\cite{Rob01}, Section 9.5 and \cite{LehCas98}, note~9.3,
pages 421--422, for general expositions and many useful references. In
particular, the Hunt--Stein theorem gives, for invariant problems,
conditions on the group (amenability) that guarantee the
existence of a minimax equivariant estimator whenever a minimax
procedure exists.~\cite{Kie57} is a key reference. All of the
examples below relate to amenable groups, such as the additive and
multiplicative groups, the group of location-scale
transformations and the group of lower triangular $p \times p$
non-singular matrices with positive diagonal elements.
\end{remark}

\section{Examples}\label{sec3}
We focus here on various applications, illustrating how the
results of Section \ref{sec2} apply to both existing and new results. We
accompany this with further observations and remarks. As
previously mentioned, such applications are quite varied with
respect to model, loss and shape of the restricted parameter
space. At the expense of some redundancy, some particular cases
are singled out (e.g., (A) is a particular case of (D), while (B)
is a particular case of (E)) for their practical or historical
importance. However, we do not focus here on specific
determinations of the MRE estimators, but do refer to textbooks
that treat in detail such topics (e.g.,
\cite{LehCas98}). Throughout, we consider loss functions that satisfy the
conditions of Proposition \ref{brown}, and our findings relate to
univariate and multivariate continuous probability models with
absolutely continuous Lebesgue densities.

\begin{enumerate}[(A)]
\item[(A)] (A single location parameter.)
Consider a location model with $X \sim f_0(x_1-\theta, \ldots,
x_n-\theta)$, with $\Omega=\Re$, known $f_0$ and invariant loss
$\rho(d-\theta)$. Consider further a~lower (or upper)
bounded parameter space $\Omega^*$ (i.e., $\Omega^*=[a,\infty)$ or
$\Omega^*=(-\infty,a]$). On one hand, $\Omega^*$ satisfies the
conditions of Theorem \ref{minimax} with the choices
$B_n=\Omega^*$, $\bar{g}_n = -n$ for $\Omega^*=[a, \infty)$ (and
$B_n=\Omega^*$, $\bar{g}_n = n$ for $\Omega^*=(-\infty,a]$). On
the other hand, following \cite{Kie57} or \cite{GirSav51} for squared error loss, the MRE or Pitman estimator
is minimax (and also Bayes with respect to the flat prior for
$\theta$ on~$\Re$). Thus Corollary \ref{mre} applies and the MRE
estimator is minimax as well for the restricted parameter space
$\Omega^*$. The result is not new (see, e.g., \cite{Kat61},
for a normal model and squared error loss; \cite{Far64}, for
strictly convex $\rho$; \cite{MarStr05N1},
for strict
bowl-shaped losses). Finally, we mention the implication that the
minimaxity property is hence shared by any dominator of the MRE
estimator, which includes quite generally the Bayes estimator of
$\mu$ associated with the flat prior on $\Omega^*$ (e.g.,
\cite{Far64,MarStr05N1}).

\item[(B)] (A single scale parameter.)
Analogously, consider scale families with densities
$\frac{1}{\sigma^n} f_1(\frac{x_1}{\sigma}, \ldots,
\frac{x_n}{\sigma})$, with natural parameter space $\Omega=
\Re^+$, known $f_1$, invariant loss $\rho(d/\sigma)$ and
restricted parameter spaces $\Omega^*=[a,\infty)$ or
$\Omega^*=(0,a]$ (with $a>0$ known). With the multiplicative
group on $\Re^+$, these restricted parameter spaces satisfy the
conditions of Theorem \ref{minimax} with $B_n=\Omega^*$ and the
choices $\bar{g}_n= \frac{1}{n}$ and \mbox{$\bar{g}_n=n$} for
$\Omega^*=[a,\infty)$ and $\Omega^*=(0,a],$ respectively. From
\cite{Kie57}, whenever a minimax estimator exists for the
unconstrained case $\sigma>0$, it is necessarily given by the MRE
estimator or equivalently by the Bayes estimator with respect to
the non-informative prior $\pi(\sigma) = \frac{1}{\sigma}
I_{(0,\infty)}(\sigma)$. Thus Theorem \ref{minimax} and
Corollary \ref{mre} apply and such MRE estimators remain minimax
for constrained parameter spaces $[a,\infty)$ and $(0,a]$; and
this quite generally with respect to model $f_1$ and loss $\rho$.

A version of the above minimaxity result for strict bowl-shaped
losses was obtained by Marchand and Strawderman \cite{MarStr05N2}. Kubokawa
\cite{Kub04} provided the result for entropy loss (i.e., $\rho
(z)=z-\log z -1$), while van Eeden \cite{van06} provided the result (actually more
general, which relates to a vector of scale parameters as in (E)
below) for scale invariant squared error loss (i.e.,
$\rho(z)=z^2$). Also, we refer to the three last references for
earlier results obtained for specific models $f_1$, namely gamma
and Fisher models. Finally, we also point out that the above
development applies to estimating powers $\sigma^r$ of $\sigma$ by
the transformation $x_i \to x_i^r$ (e.g., \cite{MarStr05N2}, for more details).

\item[(C)] (Location-scale families with the location parameter
restricted to an interval (possibly compact).)
For location-scale families with observables $X_1, \ldots, X_n$
having joint density $\frac{1}{\sigma^n}
f_2(\frac{x_1-\mu}{\sigma}, \ldots, \frac{x_n-\mu}{\sigma})$,
consider estimating $\mu$ with $\sigma>0$ (unknown) under either:
(i) the compact interval restriction $\mu\in[a,b]$, or
(ii) $\mu\in[a,\infty)$; $f_2$ known, invariant loss
$\rho(\frac{d-\mu}{\sigma})$. For (i), Theorem
\ref{minimax} applies with $B_n=\Omega^*$,
$\bar{g}_n=(-\frac{n(a+b)}{2},n)$, and $\bar{g}_n \Omega^* =
\{(\mu, \sigma) \in\Re\times\Re^+ \dvt \mu\in
[-\frac{n(b-a)}{2}, \frac{n(b-a)}{2}], \sigma
>0 \} .$ As well, Kiefer \cite{Kie57} tells us that the MRE estimator
of $\mu$ or, equivalently, Bayes with respect to the Haar right
invariant prior $\pi(\mu, \sigma) = \frac{1}{\sigma}
1_{(0,\infty)}(\sigma),$ is minimax for the unrestricted problem
with $\Omega= \Re\times\Re^+$ (subject to existence). The
conclusion derived from Corollary \ref{mre} is that
$\delta_{\mathrm{MRE}}$ is also minimax for the restricted
parameter space with $\mu\in[a,b], \sigma>0$, while a similar
development and conclusion applies for~(ii) with $B_n=\Omega^*$ and
$\bar{g}_n=(-n,1)$, a result of which also follows
from (G) below. The result for compact interval restriction
(i) generalizes the result previously obtained for scaled
squared error loss (i.e., $\rho(z)=z^2$) by Kubokawa~\cite{Kub05}.

Finally, we point out that a compact interval restriction on
$\mu$ with known $\sigma$ typically leads to a different
conclusion, with a corresponding MRE estimator that is not minimax.
A somewhat familiar justification for this (e.g., see \cite{LehCas98},
page 327 for a normal mean $\mu$ and squared error $\rho$) is as follows.
Consider $\rho$ to be strictly bowled-shaped in the sense that $\rho
'(\cdot)$ is positive
on $(0,\infty)$ and negative on $(-\infty,0)$. Denote $V_0$ and $\delta
_{\mathrm{TMRE}}$
as the constant risk of $\delta_{\mathrm{MRE}}$ and the truncation of
$\delta_{\mathrm{MRE}}$ onto the parameter
space $[a,b]$, respectively. Observe that $V_0=R(\mu, \delta_{\mathrm
{MRE}}) > R(\mu, \delta_{\mathrm{TMRE}})$ for all $\mu\in[a,b]$,
and that the compactness of the parameter space coupled with the
continuity of the risk $R(\mu, \delta_{\mathrm{TMRE}})$ imply
that $\sup_{\mu\in[a,b]} R(\mu, \delta_{\mathrm{TMRE}}) < V_0$ and
that, consequently, $\delta_{\mathrm{MRE}}$ is not minimax.

\item[(D)] (Location parameters restricted to a polyhedral cone.)
Consider independently generated copies of $X \sim
f_0(x_1-\mu_1, \ldots, x_p-\mu_p)$, with $f_0$ known, and
$\mu=(\mu_1, \ldots, \mu_p)'$ restricted to a Polyhedral cone
%
\begin{equation} \label{pcone}
\Omega_C^*=\{\mu\in\Re^p\dvt
C\mu\geq0 \} ,
\end{equation}
where $C (q \times p)$ ($q \leq p$) is of full rank (and the $0$
is a $q \times1$ vector of $0$'s). Such restricted parameter
spaces include:
\begin{enumerate}[(iii)]
\item[(i)] orthant restrictions where some or all of the
$\mu_i$'s are bounded below by $0$;
\item[(ii)] order
restrictions of the type $\mu_1 \leq\mu_2 \leq\cdots\leq\mu_r$
with $r \leq p $;
\item[(iii)] tree order restrictions of the
type $\mu_1 \leq\mu_i$ for some or all of the $\mu_i$'s;
\item[(iv)] umbrella order restrictions of the type $\mu_1 \leq\mu_2
\leq\cdots\leq\mu_m \geq\cdots\geq\mu_p$ ($m$~known).
\end{enumerate}

With $B_n=\Omega_C^*$ and $\bar{g}_n \in\Re^p$ as the additive group
elements such that $C \bar{g}_n=-n (1, \ldots, 1)'$, we obtain
$\bar{g}_n \Omega_C^*=\{\mu\in\Re^p\dvt C\mu\geq-n (1, \ldots,
1)' \}$ and choices that satisfy Theorem \ref{minimax}.
Furthermore, for invariant losses $\rho(\|d-\mu\|)$,
the results of Kiefer \cite{Kie57} tell us that, subject to risk
finiteness, the MRE or Bayes estimator for~$\mu$ with a flat prior
on $\Re^p$ is minimax for the unconstrained problem $\mu\in
\Re^p$. We infer by Corollary \ref{mre} that the same estimator is
minimax for any polyhedral cone~$\Omega_C^*$ as in (\ref{pcone}).

Other than problems in (A), the above unifies and extends several
previously established results, beginning with the Blumenthal and
Cohen \cite{BluCoh68} case of order constraints and squared error
$\rho$, and
including more recent findings by Tsukuma and Kubokawa \cite
{TsuKub08} for multivariate
normal models, the general constraint in (\ref{pcone}) and squared
error $\rho$
(also see \cite{KumSha88} and \cite
{van06}, for further
results and references). A~much-studied and important case is the
normal model with $X \sim N_p(\mu, I_p)$, \mbox{$\mu\in\Omega_C^*$}
and loss $\|d-\mu\|^2$, for which the above results apply with
\mbox{$\delta_{\mathrm{MRE}}(X)=X$}. As an interesting corollary of a result by
Hartigan \cite{Har04} and of the above, it follows that the Bayes
estimator $\delta_U$ of $\mu$ with respect to a flat prior on
$\Omega_C^*$, which Hartigan showed dominates $X$, is minimax for
$\Omega_C^*$. To conclude, we point out that a~particular case of
Hartigan's result was obtained by Blumenthal and Cohen~\cite
{BluCoh68} for
ordered location parameters in (ii) with $r=p=2$. They actually
provide a class of model densities $f_0$, including normal,
uniform and gamma densities, through conditions that ensure that
$\delta_U$ (also referred to as the Pitman estimator by
the authors) is minimax under squared error loss. They also
report on numerical evidence indicating that $\delta_U$ is not
minimax in general with respect to $f_0$.

\begin{remark}\label{rem2}
As an extension of the above, a similar development holds with the
introduction of an unknown scale parameter $\sigma$ $(\sigma>0)$.
Indeed for (at least two) independent copies from density
$\frac{1}{\sigma^p} f_2(\frac{x_1-\mu_1}{\sigma}, \ldots,
\frac{x_p-\mu_p}{\sigma})$, invariant loss
$\rho(\frac{\|d-\mu\|}{\sigma})$ and restricted parameter space
$\mu\in\Omega_C^*$, $\sigma>0$, Theorem \ref{minimax} and
Corollary \ref{mre} apply as above,\vadjust{\goodbreak} but with the MRE estimator of
$\mu$ now being generalized Bayes with respect to the prior
measure $\pi(\mu, \sigma) = \frac{1}{\sigma}
1_{(0,\infty)}(\sigma) 1_{\Re^p}(\mu)$. Moreover, if
estimating an unconstrained $\sigma$ (or $\sigma^r$) is the
objective, the MRE estimator of $ \sigma$ can be shown to be
minimax as well with the parametric restrictions (subject to risk
finiteness). This means that any minimax estimator of $\sigma^r$
for an unconstrained problem remains minimax even when $\mu\in
\Omega_C^*$.
\end{remark}

\item[(E)] (Ratios or products of scale parameters.)
For independently generated copies of $X \sim(\prod_i
\sigma_i)^{-1} f_1(\frac{x_1}{\sigma_1}, \ldots,
\frac{x_p}{\sigma_p})$ with $f_1$ a known Lebesgue density (on
($\Re^+)^p$), $\Omega=(\Re^+)^p$, consider the restriction
$\tau=\prod_i (\sigma_i)^{r_i} \geq c >0$, with the $r_i$'s known
and estimating $\tau$ under invariant loss $\rho(\frac{d}{\tau})$.
The parametric function $\tau$ includes interesting cases of
ratios $\frac{\sigma_i}{\sigma_j}$ and products $\sigma_i
\sigma_j$ (with or without nuisance parameters $\sigma_k$, $k \neq
i,j$), and the constraint on $\tau$ represents a natural scale
parameter analog of (\ref{pcone}) with $q=1$. With
$B_n=\Omega^*$, and $\bar{g}_n \in(\Re^+)^p$ the multiplicative
group element given by $\bar{g}_n=(n^{-1/r_1}, \ldots,
n^{-1/r_p})$, we obtain $\bar{g}_n \Omega^*= \{(\sigma_1, \ldots,
\sigma_p) \in(\Re^+)^p \dvt \prod_i (\sigma_i)^{r_i} \geq
\frac{c}{n^p} \}$. Thus, the conditions of Theorem~\ref{minimax}
are satisfied, and Theorem \ref{minimax} applies. Corollary
\ref{mre} applies as well, by virtue of Kiefer~\cite{Kie57}, indicating
that the MRE estimator (if it exists), or equivalently Bayes with
respect to the prior measure $\prod_i \frac{1}{\sigma_i}
1_{(0,\infty)}(\sigma_i)$, remains minimax\vspace*{1pt} for estimating $\tau$
under the lower bound constraint above. We refer to \cite{LehCas98}, Chapter 3,
problems~3.34--3.37 for examples.
Finally, with the minimax result here being quite general with
respect to the loss $\rho$ (as well as with respect to the type of
constraint and the model), we point out that the particular case
of scale-invariant squared error loss (i.e., $\rho(z)=(z-1)^2$) is
covered by van Eeden \cite{van06}, Lemma 4.5.

\item[(F)] (Linear combinations of restricted location
parameters.)
Consider location models with $X=(X_1, \ldots, X_k)' \sim\prod_i
f_i(x_i-\mu_i)$, known $f_i$'s, where we wish to estimate
$\theta\,{=}\,\sum_{i=1}^k a_i \mu_i\,{=}\,a'\mu$, under loss
$\rho(d\,{-}\,\theta)$ and the restriction
$\mu\,{\in}\,\Omega^*\,{=}\,\{\mu\,{\in}\,
\Re^n \dvt \mu_i\,{\geq}\,0\allowbreak \mbox{for }  i=1,\ldots, k \}$. For the
unconstrained version with $\mu\in\Omega=\Re^k$, the MRE
estimator (also Bayes with respect to the flat prior on $\Re^k$)
is minimax \cite{Kie57}, subject to existence and risk
finiteness. Hence, Corollary \ref{mre} (or Theorem \ref{minimax})
applies with $B_n= \Omega^*, \bar{g_n}=(-n, \ldots, -n)$
indicating that an MRE estimator remains minimax in the constrained
problem $\mu\in\Omega^*$ for estimating $\theta$. Kubokawa
\cite{Kub} has recently established the above for squared error loss,
where the MRE estimator, whenever it exists, is the unbiased
estimator $\sum_{i=1}^k a_i (X_i-b_i)$ with $E(X_i-\mu_i)=b_i$.
The result is extended here with respect to $\rho$, and
achieved with a different and more general proof.

\item[(G)] (Quantiles with parameter space restrictions.)
Consider location-scale models with $(X_1, \ldots, X_m)' \sim\frac
{1}{\sigma^m} \prod_i
f_0(\frac{x_i-\mu}{\sigma})$; $m \geq2$, $f_0$ known, $\mu\in
\Re$, $\sigma>0$; with the objective of estimating a quantile
parameter $\mu+ \eta\sigma$; (of known order
$\int_{-\infty}^{\eta} f_0(z) \,\mathrm{d}z$) under invariant loss
$\rho(\frac{d-\mu-\eta\sigma}{\sigma})$. Now,
consider restricted parameter spaces such~as:
\[
\Omega_1^*=\{(\mu, \sigma) \in\Re\times\Re^+ \dvt \mu+ \eta\sigma
\geq0 \}
\]
and
\[
\Omega_2^*=\{(\mu, \sigma) \in\Re\times
\Re^+ \dvt \mu\geq a, \sigma\geq b \geq0 \}.\vadjust{\goodbreak}
\]

Taking $B_n=\Omega_1^*$ and $\bar{g}_n=(-n,\frac{1}{n})$ such that
$\bar{g}_n \Omega_1^* = \{(\mu, \sigma) \in\Re\times\Re^+ \dvt
\mu+ \eta\sigma\geq-n \}$ and $\bar{g}_n \Omega_2^* = \{(\mu,
\sigma) \in\Re\times\Re^+ \dvt \mu\geq\frac{a}{n}- n, \sigma
\geq\frac{b}{n} \}$, we see that the conditions of Theorem~\ref{minimax} are satisfied. Moreover, subject to existence or
risk finiteness, the results of Kiefer \cite{Kie57} tell us the MRE
estimator is minimax for the unrestricted parameter space
$\Omega=\Re\times\Re^+$. Hence, Corollary \ref{mre} applies and
tells us that such MRE estimators are minimax for restricted
parameter spaces $\Omega_1^*$ and $\Omega_2^*$. Previously
studied models, for which the above results apply, include
exponential and normal $f_0$'s. For instance, consider a standard
normal $f_0$ and squared error~$\rho$, where equivariant estimators are
of the
form $\bar{X} + \eta c S$, $\delta_{\mathrm{MRE}}(X_1, \ldots, X_m) =
\bar{X} + \eta c_m
S$, with constant and minimax risk $1+
\eta^2(1-(m-1)c_m^2)$, and where $\bar{X}=\frac{1}{m}\sum_{i=1}^m X_i,
S^2=\sum_{i=1}^{m}
(X_i-\bar{X})^2$ and
$c_m=\frac{\Gamma({m}/{2})}{\sqrt{2}\Gamma({(m+1)}/{2})}$
(\cite{Fer67}, page 182).
The general result above tells us the $\delta_{\mathrm{MRE}}$
remains minimax for parameter spaces $\Omega_1^*$ and
$\Omega_2^*$ under squared error loss.

Observe also that the above development relative to $\Omega_1^*$ is
still valid whenever $\eta=0$, which relates to the problem of
estimating a median or mean for symmetric $f_0$'s, with the
corresponding minimaxity result previously obtained by Kubokawa~%
\cite{Kub04} for scale-invariant squared error loss (i.e., $\rho(z)=z^2$
above). Similar results follow with an upper bound of $0$ for $\Omega
_1^*$, as
well as an upper bound for $\mu$ and/or an upper bound for
$\sigma$ in the case of $\Omega_{2}^{*}$. Finally, we point out that
the minimaxity result and
development above follow without emendation for the general case
of non-independent components with joint density
$\frac{1}{\sigma^m} f(\frac{x_1-\mu}{\sigma}, \ldots,
\frac{x_m-\mu}{\sigma})$.

\item[(H)] (Restricted covariance matrices.)
Consider a summary statistic $S \sim\hbox{Wishart}(\Sigma,\allowbreak p,m)$
with $m \geq p$ and $\Sigma$ positive definite. Moreover, suppose
that we wish to estimate $\Sigma$ with invariant loss (under the
general linear group) $L(\Sigma, \delta) =
\psi(\Sigma^{-1} \delta)$, with $\psi(y) = \operatorname{tr}(y) - \log
|y| -p$
and $\psi(y)=\operatorname{tr}(y-I_p)^2$ as interesting examples. A standard
method to derive a minimax estimator here (e.g., \cite{Eat89},
Section 6.2) is to consider the best equivariant estimator under
the subgroup $G_{T}^+$ of lower triangular matrices with positive
diagonal elements. Such equivariant estimators can be shown to
have constant risk, be of the form $\delta_A(S)=(S^{1/2}) A
(S^{1/2})'$ with $A$ symmetric and $S^{1/2}$ the unique square
root of $S$ element belonging to $G_{T}^+$ and with the optimal
choice (MRE) being minimax. For instance, under loss
$\operatorname{tr}(\Sigma^{-1}\delta-I_p)^2$, the BEE is minimax and
given by
$\delta_{A_0}$, with $A_0$ the diagonal matrix with elements
$(m+p-2i+1)^{-1}$; $i=1,\ldots, p$; \cite{JamSte61}.

Now, consider restrictions on $\Sigma$ of the type
$\Omega^*=\{\Sigma>0 \dvt |\Sigma| \geq c_1 >0 \}$ or
$\Omega^*=\{\Sigma>0 \dvt \operatorname{tr}(\Sigma) \geq c_2 >0 \}$. It
is easy
to see in both cases that the conditions of Theorem \ref{brown}
apply with $\bar{g}_n= \frac{1}{n}I_p$ and $B_n=\Omega^*$. Hence,
the above MRE estimators remain minimax under the above
restrictions by virtue of Corollary \ref{minimax}.
\end{enumerate}

\section*{Concluding remarks}
We have provided in this paper a rich and vast collection of novel
minimax findings for restricted parameter spaces.
Furthermore, we have established a unified framework not only
applicable to many new situations, but also covering many
generalizations of existing minimax results with respect to model
and loss. For the sake of clarity and in a~summary attempt to
draw a sharper distinction between existing and new results to the
best of our knowledge, we point out or reiterate that:

\begin{itemize}
\item Results in (A) and (B) are not new except for the slight
generalization on the loss with our results here applicable to
losses that are not necessarily strictly bowl-shaped.

\item Situations (C)--(F) have been studied by
others with existing minimax results for squared error $\rho$. Our
results cover more general losses $\rho$ in all these cases.

\item Remark \ref{rem2}, situations (G) and (H), correspond for the most
part to new problems and the given results are novel.
\end{itemize}

\section*{Acknowledgements}
The authors are grateful to Larry Brown for useful discussions and
insight concerning Proposition \ref{brown}. The authors are also
grateful for constructive comments and
suggestions provided by two reviewers. Thanks also to Yogesh
Tripathi and Mohammad Jafari Jozani for stimulating exchanges
concerning our applications to linear combinations of location
parameters (F) and to quantiles (G). \'{E}ric Marchand
gratefully acknowledges NSERC of Canada, which provides financial
research support.

%

\printhistory

\end{document}